\documentclass[twoside,leqno,twocolumn]{article}
\usepackage{ltexpprt}

\usepackage{graphicx}
\usepackage{epstopdf}

\begin{document}

\title{\Large Preconditioned warm-started Newton-Krylov methods for MPC with
discontinuous control\thanks{
knyazev@merl.com,~Alexander.Malyshev@uib.no}}

\author{Andrew Knyazev\thanks{Mitsubishi Electric Research Laboratories (MERL)
201 Broadway, 8th floor, Cambridge, MA 02139, USA} \\
\and
Alexander Malyshev\thanks{University of Bergen, Department of Mathematics,
Postbox 7803, 5020 Bergen, Norway}}

\date{}

\maketitle


\begin{abstract} \small\baselineskip=9pt
We present Newton-Krylov methods for efficient numerical solution of
optimal control problems arising in model predictive control, where the
optimal control is discontinuous. As in our earlier work, preconditioned
GMRES practically results in an optimal $O(N)$ complexity, where $N$ is
a discrete horizon length. Effects of a warm-start, shifting along the
predictive horizon, are numerically investigated. The~method is tested
on a classical double integrator example of a minimum-time problem with
a known bang-bang optimal control.
\end{abstract}

\section{Introduction.}

We present and study numerically an approach to solution of minimum-time problems
by the Model Predictive Control (MPC) method. The approach has been designed in
our previous papers \cite{KnFuMa:15}-\cite{KnMa:16}. We apply it here
to the double integrator problem, which is a classical example of the minimum-time
optimal control problem. One of the main difficulties for numerical solution
of this problem is caused by that the optimal control input is a function
with large discontinuities, which makes the continuation procedures along
the trajectories harder for accurate computation.

The minimum-time optimal control problems have been well studied; see e.g.,
\cite{AtFa:06,BrHo:75,ChBaLe:12,Ga04,BrDiSW:11}.
In simple cases such as a double integrator with constraints on the control input
\cite{Lo:17}, a closed-form of the feedback law can be found analytically. In general cases,
obtaining a closed-form solution is a challenging task and a minimum-time open-loop
trajectory is generated numerically by applying either direct methods,
such as the direct shooting method, penalty function method, and SQP, or indirect
methods such as the multiple shooting method and collocation method
\cite{DiFeHa:09,Gi:13,KoFePeDi:15,Oht:04,Pa:14,RaWR:98},\cite{ShKeCo:10a}-\cite{TaOh:04},\cite{WaBo:10}.
To provide robustness to unmeasured uncertainties, the open-loop control can be augmented
with a feedback stabilizer to the computed open-loop trajectory.

Applying the MPC philosophy to the minimum-time control involves recomputing the
open-loop state and control trajectory subject to pointwise-in-time state and control
constraints, terminal state constraint, and the current state as the initial conditions.
The computed control trajectory is applied at the next time step; see, e.g.,
\cite{Ma:02,RaMa:13,ZaBi:09}.

We use a time scaling linear transformation to obtain a fixed-end time problem,
with the terminal time appearing as a multiplicative parameter in the dynamic
equations. We then convert the model to a discrete time and formulate
a prediction problem, where both the control input sequence and the
horizon length, now appearing as a parameter, have to be optimized.

We use the continuation method for advancing the predictive horizon in time.
Since the optimal control input is a piecewise continuous function and the discontinuities
are severe, the perturbations of the unknown vectors during the continuation can be
sufficiently large. There are two simple remedies for minimizing these perturbations:
(1) application of the Newton-like iterations several times and (2) the use of suitable interpolations for current data when reusing them as a warm-start in the next horizon interval. In our numerical experiments, satisfactory results are obtained by the aid of both techniques.
In some situations, the linear interpolation shows impressive improvements, in other
situations, multiple application of the Newton-like iterations is the only tool for
computations without failure.

In \cite{KnMa:16}, we propose a new preconditioner technique for the Newton-like iterations.
It works in the case, when the problem has only few state constraints, e.g., terminal state
constraints. The double integrator problem satisfies these restrictions, and we present
a construction of such preconditioner for this problem. The preconditioned iterations
have the arithmetical cost $O(N)$, where $N$ is the number of grid points
in the predictive horizon.

\section{Formulation of the prediction problem}
\label{s2}

We want to develop a numerical method for solving the double integrator problem \cite{BrHo:75},
where a bounded control input $|u|\leq1$ drives the dynamic system $\ddot{\chi}=u$
from an~initial state $\chi(t_0),\dot{\chi}(t_0)$ to a final state $\chi(t_f),\dot{\chi}(t_f)$
in minimum time $t_f-t_0$. We introduce the state variables $x_1=\chi$, $x_2=\dot{\chi}$
and convert the problem with the free terminal time $t_f$ to a problem with a fixed terminal time
by the linear time scaling $\tau = (t-t_0)(t_f-t_0)$ such that $0\leq\tau\leq1$.
The system dynamics after scaling becomes
$\frac{d}{d\tau}\left[\begin{array}{c}x_1\\x_2\end{array}\right]=
(t_f-t_0)\left[\begin{array}{c}x_2\\u\end{array}\right]$.
We also replace the inequality $|u(\tau)|\leq1$ by the equality constraint
$u^2+u_d^2-1=0$, $0\leq\tau<1$, where $u_d(\tau)$ is a dummy (or slack) argument.
Our goal is to minimize the performance index $J=p=t_f-t_0$ by choosing a suitable
control input $u(\tau)$.

A general continuous-time model for the minimum time problem after time scaling reads as follows:
given a performance index $J=\phi(\tau,x(\tau),p)|_{\tau=1}+\int_0^1L(\tau,x,u,u_d,p)d\tau$,
find $\min_{p,u}J$ subject to the equality constraints
\begin{eqnarray*}
\frac{d}{d\tau}x(\tau)=f(\tau,x(\tau),u(\tau),p),\quad 0<\tau<1,\\
x(0)=x_0,\quad \psi(\tau,x(\tau),p)|_{\tau=1}=0,\\
C(\tau,x(\tau),u(\tau),u_d(\tau),p)=0,\quad 0\leq\tau\leq1,
\end{eqnarray*}
where $x\in R^{n_x}$ and $u\in R^{n_u}$ are the state and control input, respectively.
The vector function $C$ determines equality constraints, where possible dummy arguments
$u_d$ correspond to inequality constraints. For instance, an inequality constraint
$c_j(t,x,u)\leq0$ can be replaced by the equality constraint $c_j(t,x,u)+u_d^2=0$ or by $c_j(t,x,u)+exp(u_d)=0$, and then a suitable penalty term, e.g., $-w_dp\int_0^1u_dd\tau$ must be added to the performance index $J$, where $w_d>0$ is a small constant required by the interior point method \cite{WaBo:10}.

Let us consider a discretization of the continuous-time model on the uniform
grid $\tau_i=i\Delta\tau$, where $\Delta\tau=1/N$. The adjustable variables are the control time sequence $u(0)$, $u(1)$, \ldots, $u(N-1)$ and the parameter $p$ that need to be simultaneously optimized. The discrete-time optimal control problem can be recast in the following form:
given a performance index
\[
J = \phi(\tau_N,x_N,p) + \sum_{i=0}^{N-1}L(\tau_i,x_i,u_i,u_{d,i},p)\Delta\tau,
\]
find
\[
\min_{u_i,p} J,
\]
subject to a fixed initial value $x_0$ and the equality constraints
\[
x_{i+1} = x_i + f(\tau_i,x_i,u_i,p)\Delta\tau,\quad i = 0,1,\ldots,N-1,
\]
\[
C(\tau_i,x_i,u_i,u_{d,i},p) = 0,\quad  i = 0,1,\ldots,N-1,
\]
\[
\psi(\tau_N,x_N,p) = 0.
\]
From several choices for the functions $\phi$, $L$ and $\psi$ in the double integrator problem.
we have tried two: (1) $\phi(\tau_N,x_N,p)=p$, $\psi(\tau_N,x_N,p)=x(\tau_N)-x_f$,
$L(\tau_i,x_i,u_i,u_{d,i},p)=-w_du_{d,i}p$, and (2) $\phi(\tau_N,x_N,p)=p+\frac{\beta}{2}x_N^2$,
$L(\tau_i,x_i,u_i,u_{d,i},p)=-w_du_{d,i}p+\alpha u_i^2p$ for a small $\alpha>0$ and large $\beta>0$.

To solve the discrete-time optimal control problem by the classical Lagrange multiplier
method, we introduce the costate $\lambda$ for the dynamics equations,
the Lagrange multiplier $\mu$ for the equality constraints, and Lagrange multiplier $\nu$
for the terminal constraint $\psi$ and derive necessary optimality conditions are
by means of the discrete Lagrangian
\begin{eqnarray*}
\lefteqn{\mathcal{L}(X,U)=\phi(\tau_N,x_N,p)+\sum_{i=0}^{N-1}L(\tau_i,x_i,u_i,u_{d,i},p)\Delta\tau}\\
&&{}+\lambda_0^T[x_t-x_0]\\
&&+\sum_{i=0}^{N-1}\lambda_{i+1}^T[x_i-x_{i+1}
+f(\tau_i,x_i,u_i,p)\Delta\tau]\\
&&+\sum_{i=0}^{N-1}\mu_i^TC(\tau_i,x_i,u_i,u_{d,i},p)\Delta\tau+\nu^T\psi(\tau_N,x_N,p),
\end{eqnarray*}
where
\[
X = [x_0^T,\ldots,x_{N}^T,\lambda_0^T,\ldots,\lambda_{N}^T]^T
\]
and
\begin{eqnarray*}
U = [u_0^T,u_{d,0}^T,\mu_0^T\ldots,u_i^T,u_{d,i}^T,\mu_i^T,\ldots,\\
u_{N-1}^T,u_{d,N-1}^T,\mu_{N-1}^T,\nu^T,p^T]^T.
\end{eqnarray*}
The optimality conditions are given by the stationarity conditions
\[
\frac{\partial \mathcal{L}^T}{\partial X}(X,U)=0 \mbox{ and }
\frac{\partial \mathcal{L}^T}{\partial U}(X,U)=0.
\]
Further derivations use the Hamiltonian function
\begin{eqnarray*}
\lefteqn{H(t,x,\lambda,u,u_d,\mu,p) = L(t,x,u,u_d,p)}\hspace*{5em}\\
&&{}+\lambda^Tf(t,x,u,p)+\mu^TC(t,x,u,u_d,p).
\end{eqnarray*}

The optimality conditions are reformulated below as a system of nonlinear equations
\begin{equation}\label{e1}
 F[U(t),x_t,t]=0
\end{equation}
with respect to the unknown vector $U(t)$. We recall that in the MPC method $x_t$
denotes the current measured state vector, which serves as the initial value $x_0$
in the prediction problem, and $t$ is the current system time.
The vector function $F(U,x,t)$ is obtained from the optimality conditions after elimination
of all states $x_i$ and costates $\lambda_i$ by the following procedure.
\begin{enumerate}
\item Having the current measured state $x_t$, set $x_0=x_t$ and compute
$x_i$, $i=1,2\ldots,N$, by the forward recursion
\[
x_{i+1} = x_i + f(\tau_i,x_i,u_i,p)\Delta\tau,\, i=0,\ldots,N-1.
\]
Then starting with the value
\[
\lambda_N=\frac{\partial\phi^T}{\partial x}(\tau_N,x_N,p)+
 \frac{\partial\psi^T}{\partial x}(\tau_N,x_N,p)\nu
\]
compute the costate $\lambda_i$, $i=N\!-\!1,\ldots,0$, by the backward recursion
\[
\lambda_i=\lambda_{i+1}+\frac{\partial H^T}{\partial x}
(\tau_i,x_i,\lambda_{i+1},u_i,u_{d,i},\mu_i,p)\Delta\tau.
\]
\item Using the values $x_i$ and $\lambda_i$, $i=0,1\ldots,N$, computed
 by the forward and backward recursions, compute the vector function $F[U,x,t]$
 as the vector $\frac{\partial \mathcal{L}^T}{\partial U}(X,U)$. Explicit
 formulas for the components of $F[U,x,t]$ can be found in \cite{KnMa:16}.
\end{enumerate}

\section{The Newton-Krylov method}

The system of nonlinear equations (\ref{e1}) with respect to $U(t)$ must be solved in real time
at each instance $t$ of the discrete system time. Assume that $U^{(0)}(t)$ is an approximation
to the exact solution $U(t)$. Then $F[U(t),x(t),t]-F[U^{(0)}(t),x(t),t]=-F[U^{(0)}(t),x(t),t]$.
Let us denote $b=-F[U^{(0)}(t),x(t),t]$ and $\Delta U=U(t)-U^{(0)}(t)$.
For a sufficiently small scalar $h>0$, e.g. $h=10^{-8}$, which may be different from the system time step $\Delta t$ and from the horizon time step~$\Delta \tau$, we introduce the finite difference operator
\[
a(V)=(F[U^{(0)}(t)+hV,x(t),t]-F[U^{(0)}(t),x(t),t]).
\]
Equation (\ref{e1}) is then equivalent to the operator equation
\[
h a(\Delta U/h)=b, \mbox{ where } b=-F[U^{(0)}(t),x(t),t].
\]
We suppose that $h$ is sufficiently small such that the operator $a(V)$ is almost linear in the
sense that $\|h a(\Delta U/h)-a(\Delta U)\|\leq\delta\|\Delta U\|$ with a tiny constant $\delta>0$.
Thus, solution of (\ref{e1}) is reduced to solving the operator equation
\[
a(\Delta U)=b \mbox{ where } b=-F[U^{(0)}(t),x(t),t]
\]
for a given approximation $U^{(0)}(t)$ and setting $U(t)=U^{(0)}(t)+\Delta U$.

Given formulas for computing the vector function $F[U,x_t,t]$, we want to solver
equation (\ref{e1}) at the points of the uniform grid $t_j=t_0+j\Delta t$,
$i=0,1,\ldots$ by using the above proposed approach.

At the initial state $x_0=x(t_0)$, we find an approximate  solution $U_0$ to the equation $F[U_0,x_0,t_0]=0$ by a suitable optimization procedure. The dimension of the vector $u(t)$
is denoted by $n_u$. Since
\[
U(t)=[u_0^T,u_{d,0}^T,\mu_0^T,\ldots,\ldots,\nu^T,p^T]^T,
\]
the first block entry of $U_0$, formed from the first $n_u$ elements of $U_0$,
is taken as the control $u_0$ at the state $x_0$.
The next state $x_1=x(t_1)$ is either measured by a sensor or computed
by the formula $x_1=x_0+\Delta tf(t_0,x_0,u_0)$.

The computation of $U_j$ for $j>0$ is as follows.
By $e_k$, we denote the $k$-th column of the $m\times m$
identity matrix, where $m$ is the dimension of the vector $U_j$.
At time $t_j$, we arrive with the state $x_j$ and the vector $U_{j-1}$.
The vector $U_{j-1}$ serves as the approximation $U^{(0)}_j$, and such
an approximation is often called the warm start.
The difference operator
\[
a_j(V)=\left(F[U_{j-1}+hV,x_j,t_j]-F[U_{j-1},x_j,t_j]\right)/h,
\]
implicitly determines an $m\times m$ matrix $A_j$ with the columns
\[
A_je_k=a_j(e_k),\, k=1,\ldots,m,
\]
which approximates the symmetric Jacobian matrix $F_U[U_{j-1},x_j,t_j]$
so that $a_j(V)=A_jV+O(h)$.
At the current time $t_j$, our goal is to solve the equation
\begin{equation}\label{e2}
 ~~~~a_j(\Delta U_j)=b_j, \mbox{ where } b_j=-F[U_{j-1},x_j,t_j].
\end{equation}
Then we set $U_j=U_{j-1}+\Delta U_j$ and choose the first $n_u$ components of $U_j$ as
the control $u_j$. The next state $x_{j+1}=x(t_{j+1})$ either comes from a sensor, estimated,
or computed by the formula $x_{j+1}=x_j+\Delta tf(t_j,x_j,u_j)$.

Equation (\ref{e2}) can be solved approximately by generating the matrix $A_j$
and then solving the system of linear equations $A_j\Delta U_j=b_j$ using a
direct method, e.g.,\ the Gaussian elimination.

An alternative method of solving (\ref{e2}) is by a Krylov subspace iteration,
e.g.,\ by GMRES method \cite{Sa:03}, for which we do not generate the matrix~$A_j$ explicitly.
Namely, we simply use the operator $a_j(V)$ instead of computing the matrix-vector product $A_jV$,
for arbitrary vectors $V$; cf.,~\cite{Ke:95}. The arithmetic cost of the alternative method
can be often less than that of the direct method because the matrix $A_j$ is not generated and
the $O(N^3)$ cost due to the Gaussian elimination is avoided.

We refer to the Newton-like refinement by solving equation (\ref{e2}) with the GMRES method
as the Newton-Krylov method.

\paragraph{Remark 1.}
The above numerical procedure uses the so-called warm start. Specifically, the operator
$a_j(V)$ and the right-hand-side $b_j$ are computed at time $t_j$ using the value $U_{j-1}$ from
the previous time $t_{j-1}$ as an approximation $U^{(0)}_j$.
This approximation does not take advantage of the structure of the underlying problem.
In our case, the vector $U(t)$ contains the components $u(\tau_i)$, $u_d(\tau_i)$,
$\mu(\tau_i)$, $\nu_1$, $\nu_2$ and $p$ computed in the prediction horizon $[0,1]$.
It is possible to apply, e.g., the linear interpolation of the variables $u(\tau_i)$, $u_d(\tau_i)$,
$\mu(\tau_i)$ from the previous horizon interval $[t_{j-1},t_{j-1}+p_{j-1}]$ to the current
horizon interval $[t_j,t_j+p_j]$, where $p_{j-1}$ equals $p$ from $U_{j-1}$ and $p_j=p_{j-1}-\Delta t$.
Such an interpolation is referred to as the shifting along the predictive horizon;
see a related discussion in \cite{DiFeHa:09}.

\paragraph{Remark 2.}
The shifting along the predictive horizon, described in Remark 1, is a light improvement,
which may help in some situations. A more cardinal improvement of the warm start
is solving the operator equation $a(\Delta U)=b$ and subsequent correction $U+\Delta U$
several times, say, 2-3 Newton-like iterations instead of one.

\section{Sparse preconditioner}
\label{sec:sprec}

Convergence of the GMRES iteration can be rather slow if the matrix $A_j$ is ill-conditioned.
However, the convergence can be improved by preconditioning; see e.g. \cite{Sa:03}. A matrix $T$ that approximates a matrix $A^{-1}$ and such that computing the product $Tr$ for an arbitrary vector $r$ is relatively easy, is referred to as a preconditioner. The preconditioning for the system of linear
equations $Ax=b$ formally replaces the original system $Ax=b$ with the equivalent preconditioned linear system $TAx=Tb$. If the condition number $\kappa(TA)=\|TA\|\|A^{-1}T^{-1}\|$ of the matrix $TA$ is small, convergence of iterative solvers for the preconditioned system become fast.
Instead of $T$, its inverse $M=T^{-1}$ is often called a preconditioner.

Our specific optimization problem over a predictive horizon admits construction of
efficient preconditioners $M_j$ with a sparse structure. These preconditioners
have been introduced in \cite{KnMa:16}, and we refer to this paper for a more detailed description.

The construction of $M_j$ is based on several observations.
The first one is the fact that the Jacobi matrix $F_U$ is the Schur complement of the Jacobi matrix for the Lagrangian $\mathcal{L}$. If we denote the result of the forward and backward recursions by $X=g(U)$, then most entries of $F_U$ are approximated by the matrix $\mathcal{L}_{UU}(g(U),U)$ due to the second fact, which affirms that the states $x_i$, computed by the forward recursion, and
the costates $\lambda_i$, computed by the subsequent backward recursion, satisfy
the following property: $\partial x_{i_1}/\partial u_{i_2}=O(\Delta\tau)$,
$\partial\lambda_{i_1}/\partial u_{i_2}=O(\Delta\tau)$, $\partial x_{i_1}/\partial\mu_{i_2}=0$ and
$\partial\lambda_{i_1}/\partial\mu_{i_2}=O(\Delta\tau)$.
The second fact is a corollary of theorems about the derivatives of solutions of ordinary
differential equations with respect to a parameter; see, e.g.,~\cite{Pon:62,Fi:88}.
The matrix $\mathcal{L}_{UU}(g(U),U)$ is very sparse and easily computable.
We have to generate explicitly by the formula $A(e_k)$ only the last $l$ columns and rows of $F_U$,
where the integer $l$ denotes the sum of dimensions of $\psi$ and $p$.

As a result, the setup of $M_j$, computation of its LU factorization, and application
of the preconditioner all cost $O(N)$ floating point operations. The memory requirements are also of order $O(N)$.

\section{Simulations of a Double Integrator System}

To evaluate the effectiveness of the Newton-Krylov method, we compute a solution to the double integrator system with control input constraints: $\min_{u(t)}T$ subject to $\ddot{\chi}=u$, $\chi(0)=x_0$, $\dot{\chi}(0)=y_0$, $\chi(T)=x_f$, $\dot{\chi}(T)=y_f$, $|u|\leq 1$.
We consider two discrete-time models of this problem and solve both numerically by the MPC approach,
where the predictive horizon is the interval $[t,T]$, $t$ being the current system time.

\textbf{Model 1} of the predictive problem on the scaled horizon has the form
\[
\min_{u,p}p\left[1- w_d\sum_{i=0}^{N-1}\Delta\tau u_d(i)\right],
\]
subject to
\[
\left[\begin{array}{c}x(i+1)\\y(i+1)\end{array}\right] =
\left[\begin{array}{c}x(i)\\y(i)\end{array}\right] +
p\Delta\tau\left[\begin{array}{c}y(i)\\u(i)\end{array}\right],
\]
\[
x(0)=x_t,\quad y(0)=y_t,\quad x(N)=x_f,\quad y(N)=y_f,
\]
\[
u^2(i)+u_d^2(i)=1.
\]

\textbf{Model 2} of the predictive problem on the scaled horizon has the form
\begin{eqnarray*}
\lefteqn{\min_{u,p}\frac{\alpha_1}{2}(x(N)-x_f)^T(x(N)-x_f)}\\
&&{}+p\left[1- w_d\sum_{i=0}^{N-1}\Delta\tau u_d(i)\right]
+\frac{\alpha_2}{2}p\sum_{i=0}^{N-1}\Delta\tau u(i)^Tu(i),
\end{eqnarray*}
subject to
\[
\left[\begin{array}{c}x(i+1)\\y(i+1)\end{array}\right] =
\left[\begin{array}{c}x(i)\\y(i)\end{array}\right] +
p\Delta\tau\left[\begin{array}{c}y(i)\\u(i)\end{array}\right],
\]
\[
x(0)=x_t,\quad y(0)=y_t,
\]
\[
u^2(i)+u_d^2(i)=1.
\]

We choose $\alpha_1=10^3$ and $\alpha_2=0.1$ in Model 2. The initial state is $x_0=-1$, $y_0=0$,
and the terminal state is $x_f=0$, $y_f=0$ in both models.
The interior point penalty parameter $w_d=0.005$ and is fixed as suggested, e.g.,
in paper \cite{WaBo:10}.

Below we give detailed formulas for an implementation of Model 1.
Similar formulas can be easily derived for Model 2.

The components of the discretized problem on the predictive scaled horizon $[0,1]$
are as follows:
\begin{itemize}
\item $\Delta\tau=1/N$, $\tau_i=i\Delta\tau$;
\item the unknown variables are the state $\left[\begin{array}{c}
x_i\\y_i\end{array}\right]$, the costate $\left[\begin{array}{c}
\lambda_{1,i}\\\lambda_{2,i}\end{array}\right]$, the control $u_i$,
the dummy variable $u_{d,i}$, the Lagrange multipliers
$\mu_i$ and $\left[\begin{array}{c}\nu_{1}\\\nu_{2}\end{array}\right]$,
the parameter $p$;
\item the state is determined by the forward recursion
\[
\left\{\begin{array}{l} x_{i+1}=x_i+\Delta\tau py_{i},\\
y_{i+1}=y_i+\Delta\tau pu_{i},\end{array}\right.
\]
where $i=0,1,\ldots,N-1$;
\item the costate is determined by the backward recursion ($\lambda_{1,N}=\nu_1$,
$\lambda_{2,N}=\nu_2$)
\[
\left\{\begin{array}{l} \lambda_{1,i}=\lambda_{1,i+1},\\
\lambda_{2,i} = \lambda_{2,i+1}+\Delta\tau p\lambda_{1,i+1},\end{array}\right.
\]
where $i=N-1,N-2,\ldots,0$;
\item the system of equations $F(U,x_t,y_t,t)=0$, where
\[
\hspace{-1em}U=[u_0,u_{d,0},\mu_0,\ldots,u_{N-1},u_{d,N-1},\mu_{N-1},
\nu_1,\nu_2,p]
\]
has the following rows from the top to bottom:
\[
\hspace{-3em}\left\{\begin{array}{l}
\Delta\tau\left[p\lambda_{2,i+1}+2u_i\mu_i\right] = 0 \\[1ex]
\Delta\tau\left[2\mu_iu_{d,i}-w_{d}p\right] = 0 \\[1ex]
\Delta\tau\left[u_i^{2}+u_{d,i}^2-1\right]=0
\end{array}\right.(i=0,\ldots,N-1)
\]
\[
\hspace{-1em}\left\{\begin{array}{l} x_N-x_f=0\\y_N-y_f=0\end{array}\right.\hspace{15em}
\]
\[
\hspace{-1em}\left\{\begin{array}{l}\Delta\tau\left[\sum\limits^{N-1}_{i=0}
(y_i\lambda_{1,i+1}+u_i\lambda_{2,i+1})-w_du_{d,i}\right]+1 = 0.\end{array}\right.
\]
\end{itemize}

The sparse preconditioner $M_j$ is a symmetric matrix with the sparsity structure
\[
M_j = \left[\begin{array}{ccccc}M_{11}&M_{12}\\M_{21}&M_{22}\end{array}\right],
\]
where $M_{11}$ is a block diagonal matrix with the nonsingular diagonal
$3\times3$-blocks, $i=0,\ldots,N-1$,
\[
2\Delta\tau\left[\begin{array}{ccc}
\mu_i&0&u_i\\0&\mu_i&u_{d,i}\\u_i&u_{d,i}&0
\end{array}\right],
\]
The matrix  block $\left[\begin{array}{c}M_{12}\\M_{22}\end{array}\right]$
consists of 3 columns and is computed by means of the operator $a_j(V)$ as
\[
\left[a_j(e_{3N+1}),a_j(e_{3N+2}),a_j(e_{3N+3})\right],
\]
where $e_{3N+1},e_{3N+2},e_{3N+3}$ are the columns $3N+1$, $3N+2$ and $3N+3$ of the unit matrix
of order $3N+3$. The block $M_{21}$ is transpose to $M_{12}$, $M_{21}=M_{12}^T$.

The LU factorization of $M_j$ can be computed with $O(N)$ floating point operations,
e.g., as
\[
M_j=\left[\begin{array}{ccccc}I&0\\M_{21}M_{11}^{-1}&I\end{array}\right]
\left[\begin{array}{ccccc}M_{11}&M_{12}\\0&S_{22}\end{array}\right],
\]
where $S_{22}=M_{22}-M_{21}M_{11}^{-1}M_{12}$. The setup and application
of the preconditioner also cost $O(N)$ operations.

\section{Numerical results}

All the plots are computed by the GMRES method with the absolute tolerance $10^{-6}$.

\begin{figure}
\center
\includegraphics[width=20em]{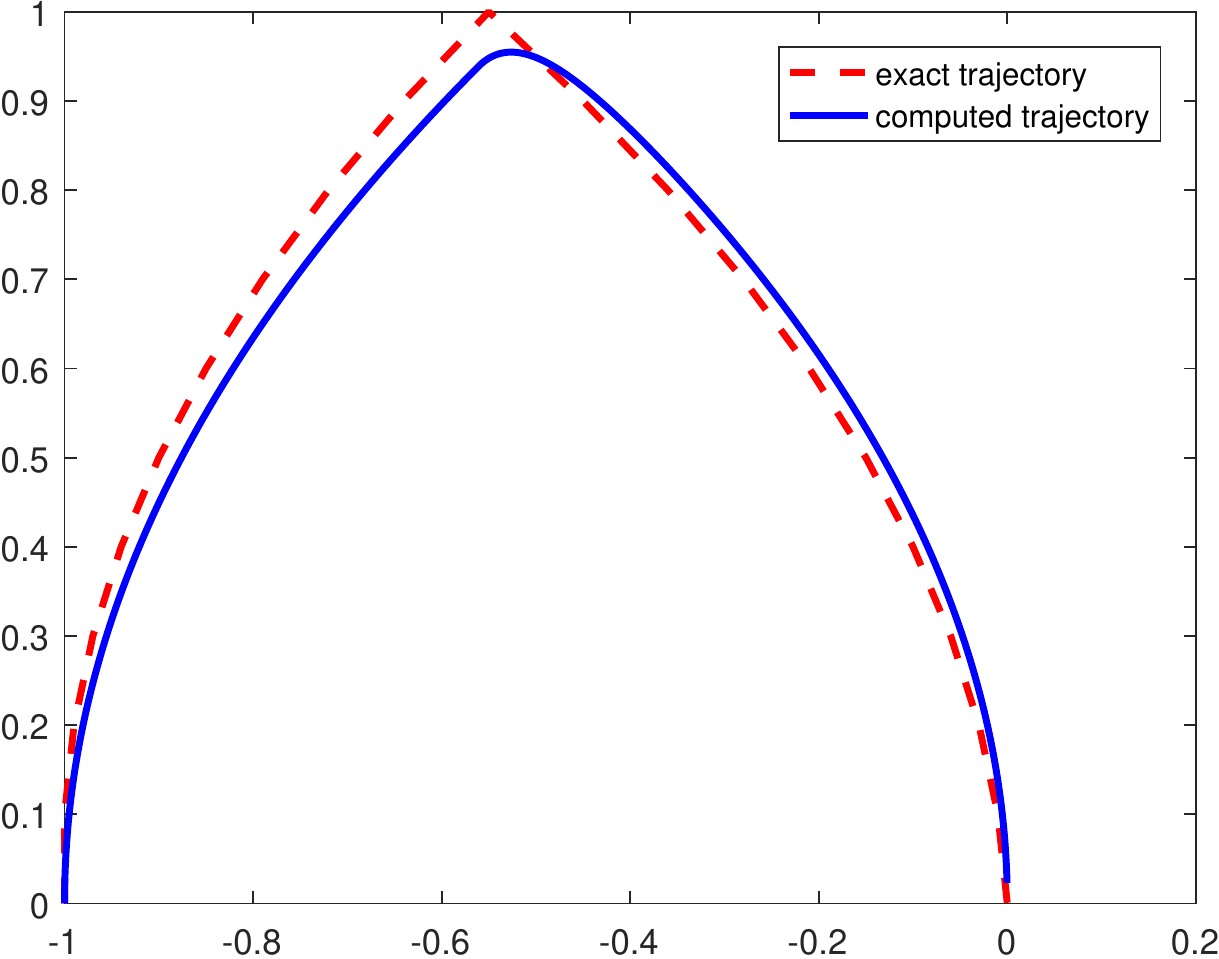}
\caption{Model 1: nominal (dashed line) vs. computed (solid line) trajectories.}
\label{fig1}
\end{figure}

\begin{figure}
\center
\includegraphics[width=20em]{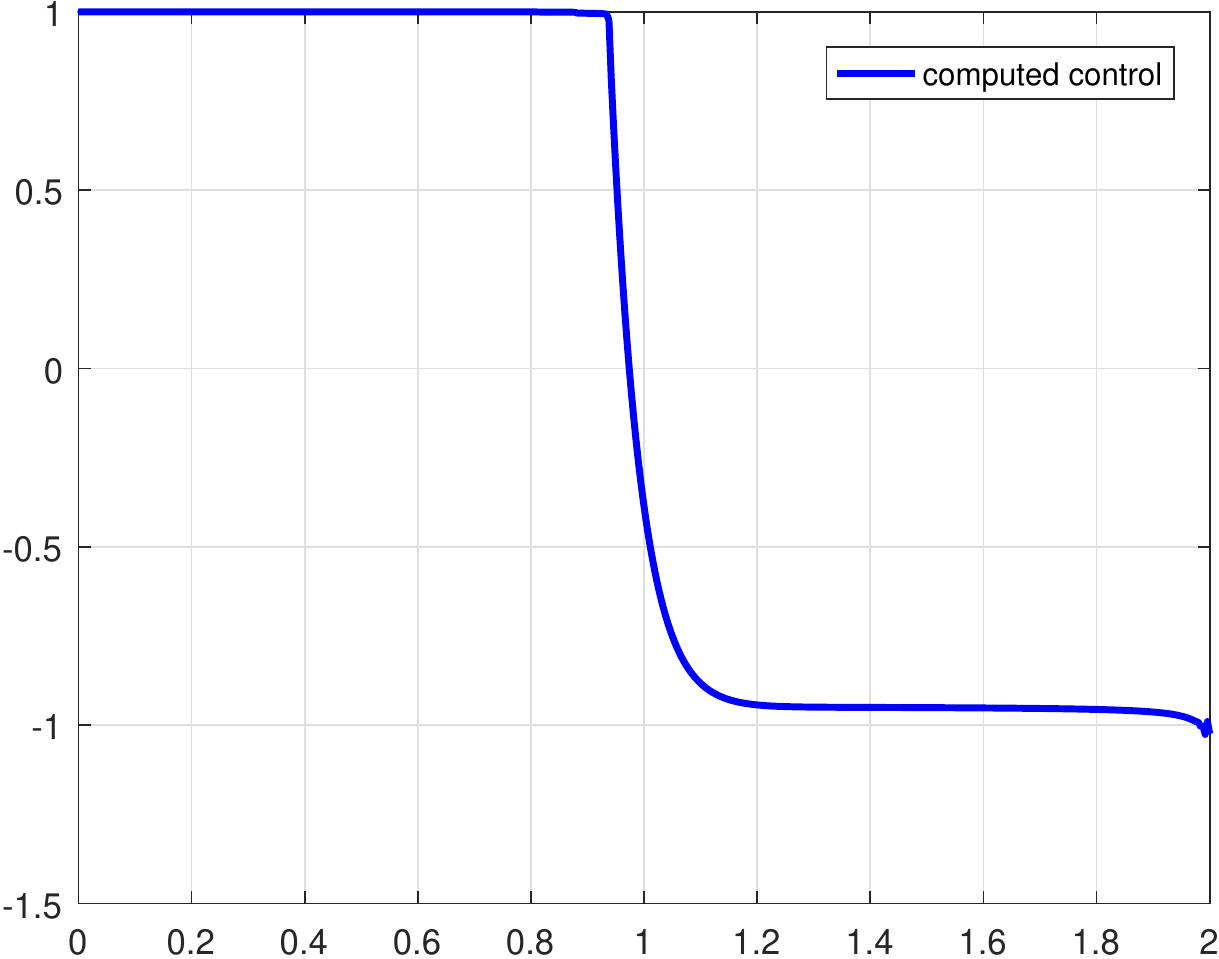}
\caption{Model 1: computed control.}
\label{fig2}
\end{figure}

\begin{figure}
\center
\includegraphics[width=20em]{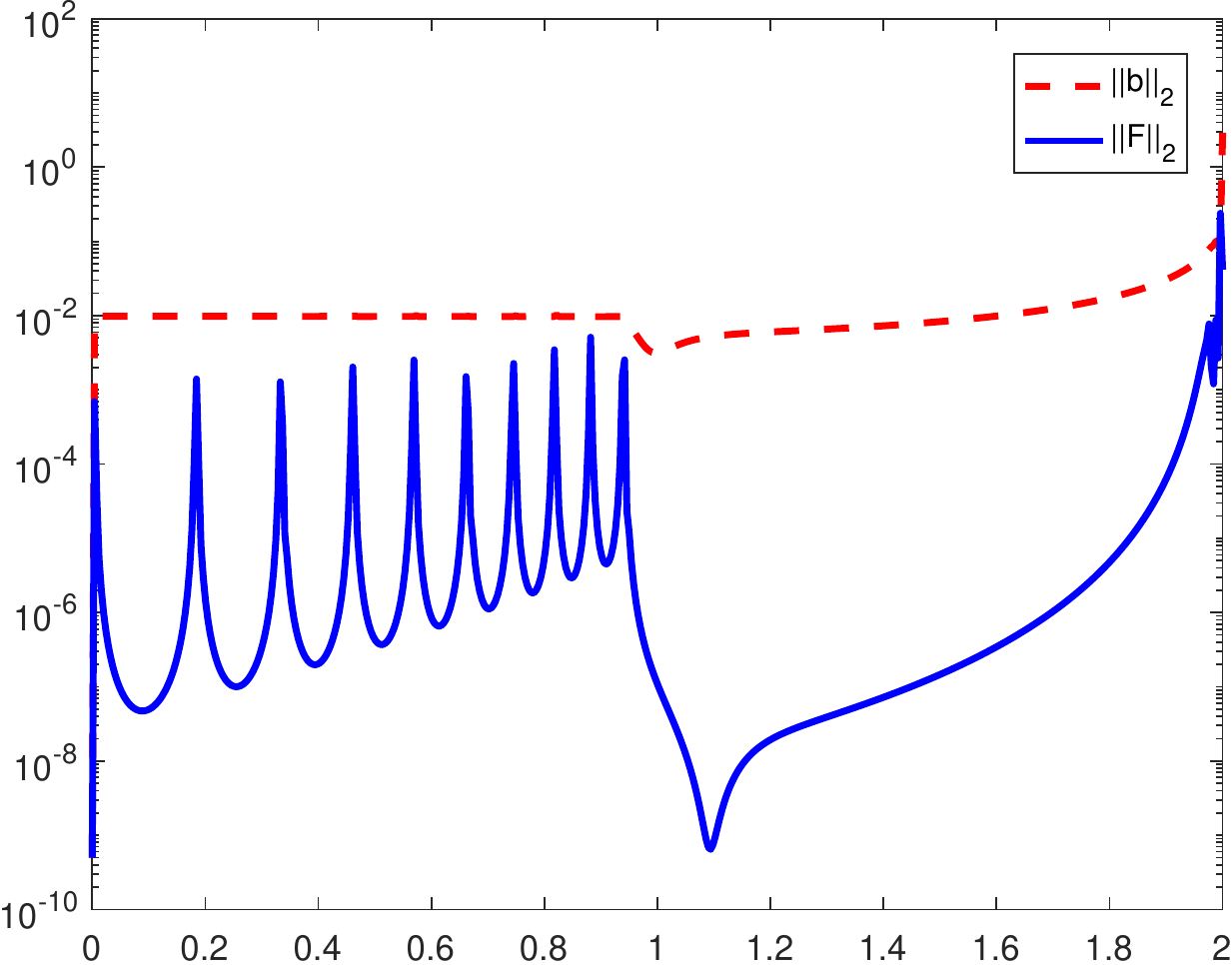}
\caption{Model 1: 2-norm of $F(U_j,x_j,t_j)$ before (dashed line) and
after (solid line) Newton's refinement.}
\label{fig3}
\end{figure}

\begin{figure}
\center
\includegraphics[width=20em]{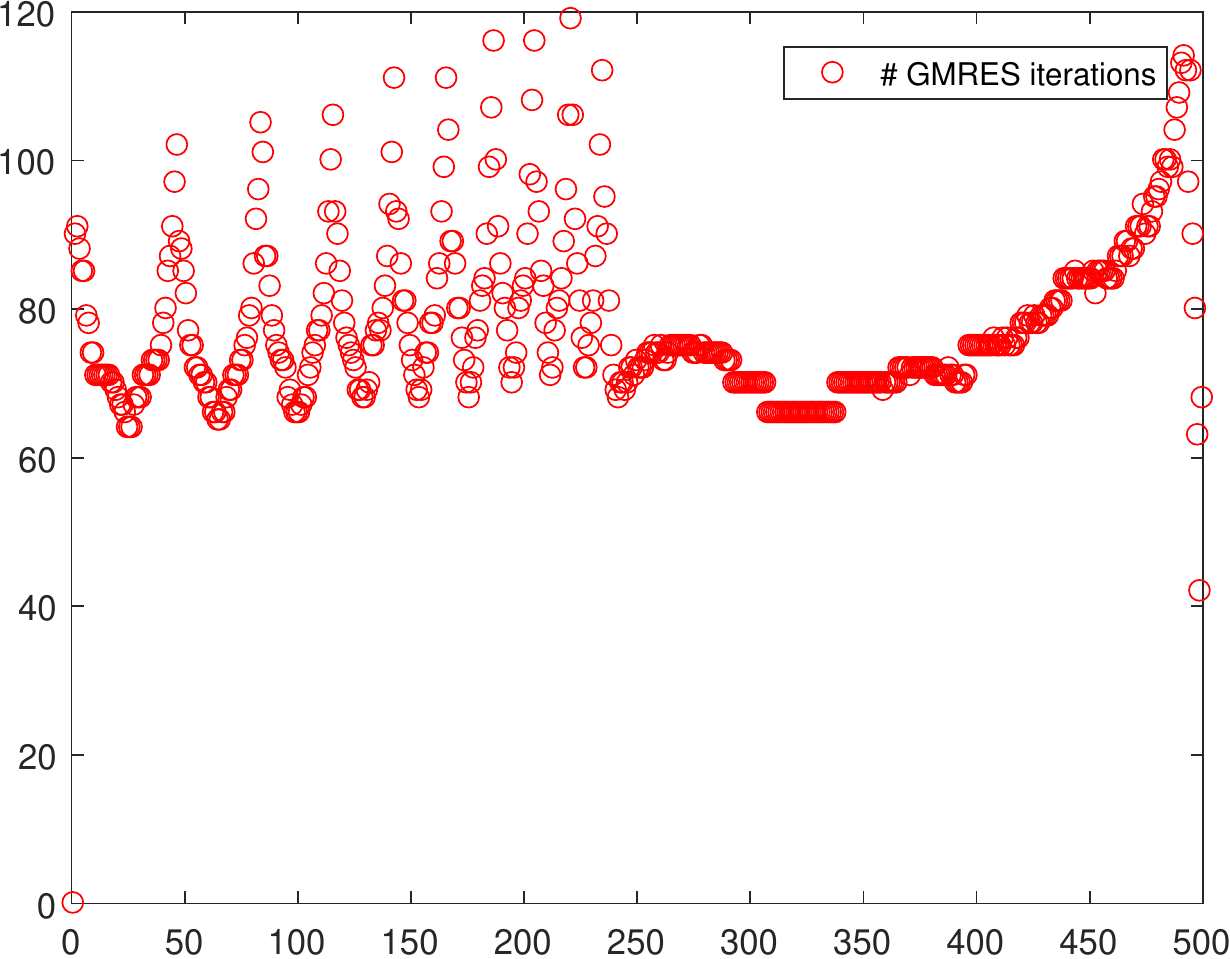}
\caption{Model 1: number of GMRES iterations without preconditioning
(77 iterations per step).}
\label{fig4}
\end{figure}

Figures \ref{fig1}, \ref{fig2}, \ref{fig3}, \ref{fig4} show the numerical results for Model 1,
where the horizon length $N$ equals 20. The exact and computed time to destination is $t_f=2$.
The system time between $t_0=0$ and $t_f$ is uniformly partitioned into 500 subintervals
with the time step $t_{j+1}-t_j=2/500$. Our MATLAB implementation fails when only 1 Newton's
refinement is used at each time instance $t_j$, and it works with 2 Newton's refinements per step.
The failure is not cured by the shifting along the predictive horizon described in Remark 1.
When the interval $[0,t_f]$ is partitioned into 1000 subintervals, our MATLAB program works with 1 Newton's refinement at each step $t_j$.

Figure \ref{fig3} shows the Euclidean norm of the residual $F[U_{j-1},x_j,t_j]$ denoted by $\|b\|_2$
and the Euclidean norm $\|F\|_2$ of the residual $F[U_{j},x_j,t_j]$ after all Newton's refinemens
at time $t_j$ (2 refinements in Figure~\ref{fig3}).

Figure \ref{fig4} displays the total number of GMRES iterations at each time $t_j$
when the GMRES method is applied without preconditioner. The average number of GMRES iterations
per step is about 77. The GMRES method with preconditioning uses only 2 iterations
per step thus demonstrating more than the 35x speedup. We recall that our preconditioning
provides the $O(N)$ complexity of the prediction problem.

The interpolation by shifting along the predictive horizon from Remark~1 does not give
an essential improvement in accuracy for Model 1 with the chosen time steps in the horizon
and the system time. The multiple application of Newton's refinements, however, provides an
essential improvement in accuracy and computation without failure.

\begin{figure}
\center
\includegraphics[width=20em]{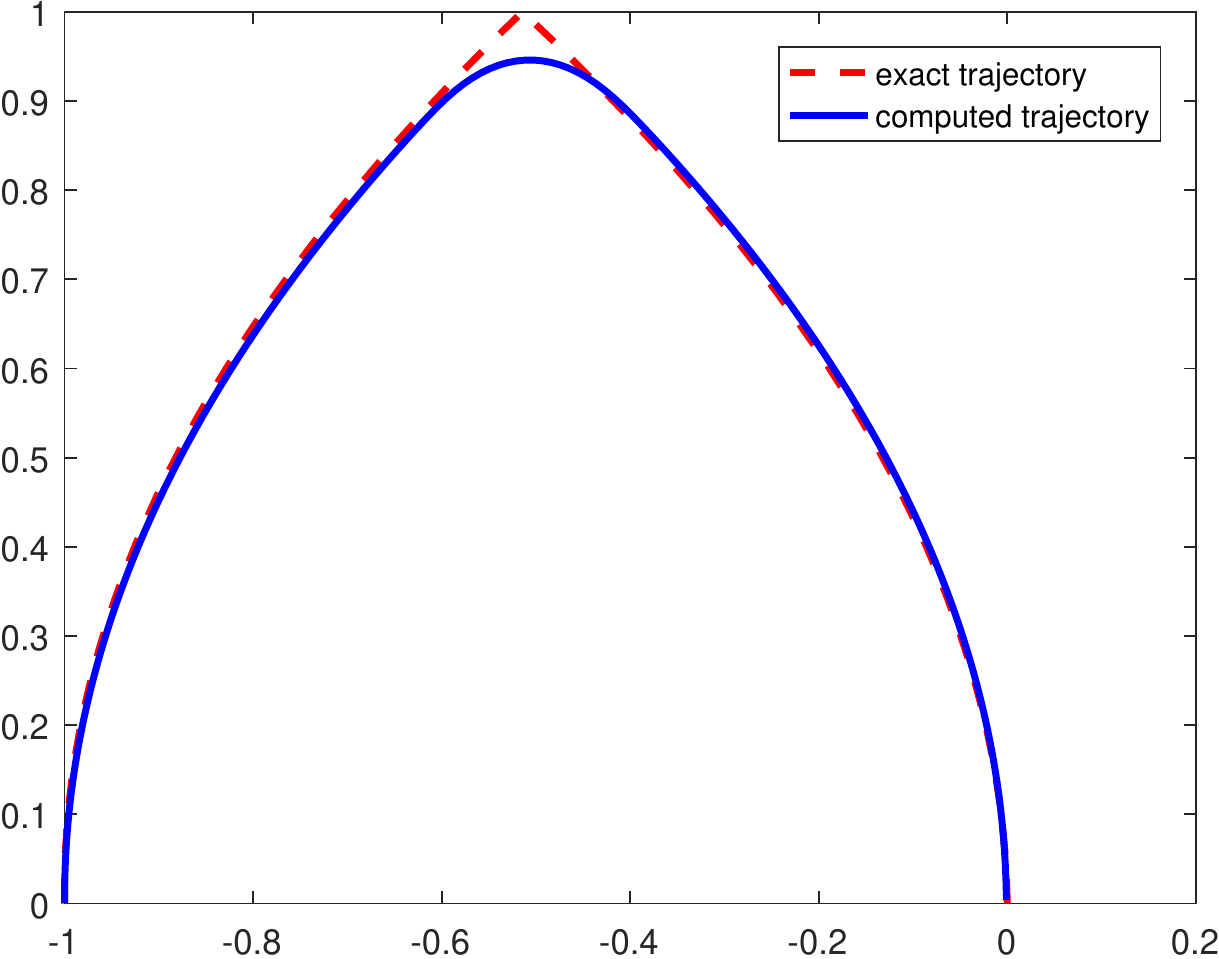}
\caption{Model 2: nominal (dashed line) vs. computed (solid line) trajectories.}
\label{fig5}
\end{figure}

\begin{figure}
\center
\includegraphics[width=20em]{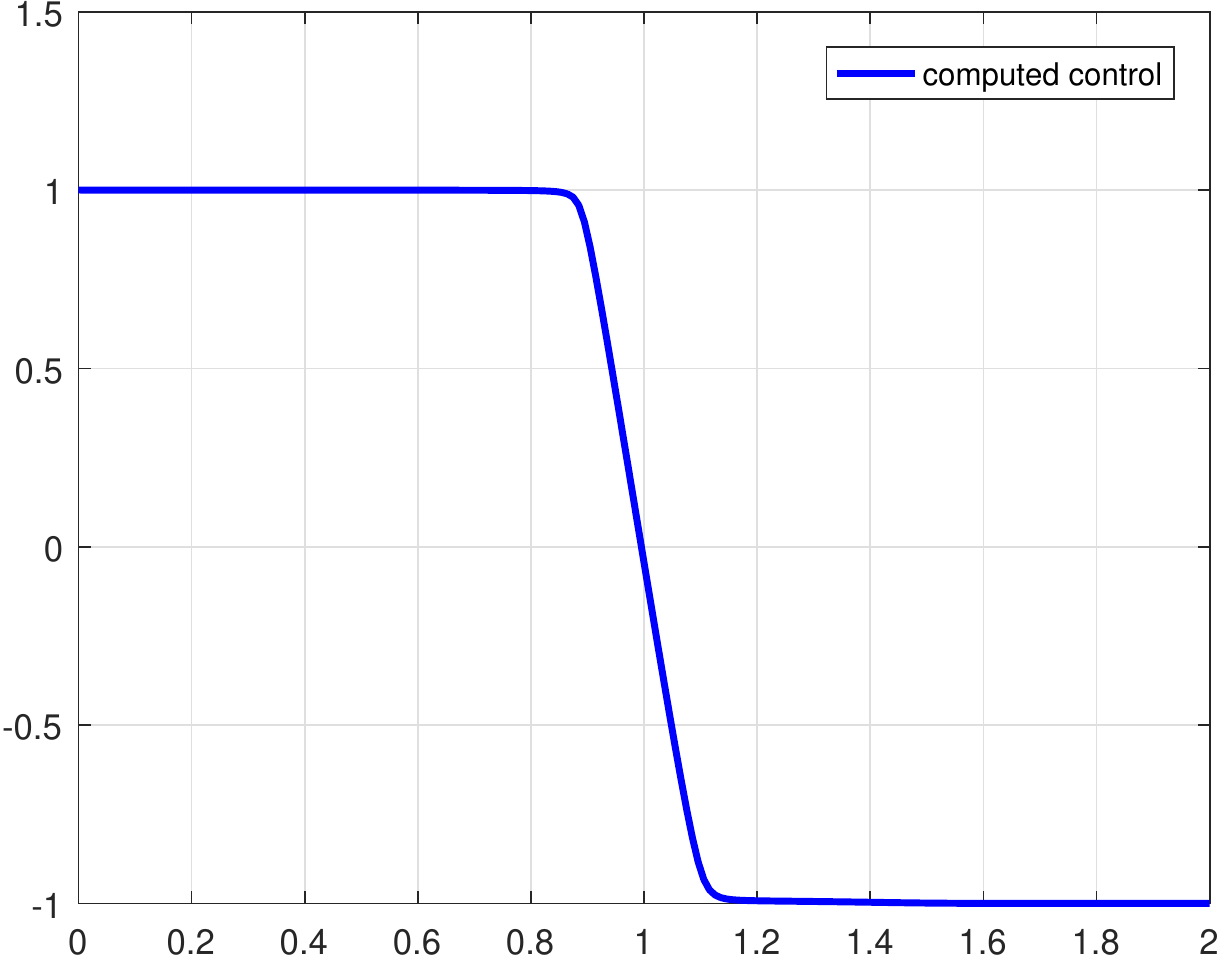}
\caption{Model 2: computed control, $N=70$.}
\label{fig6}
\end{figure}

\begin{figure}
\center
\includegraphics[width=20em]{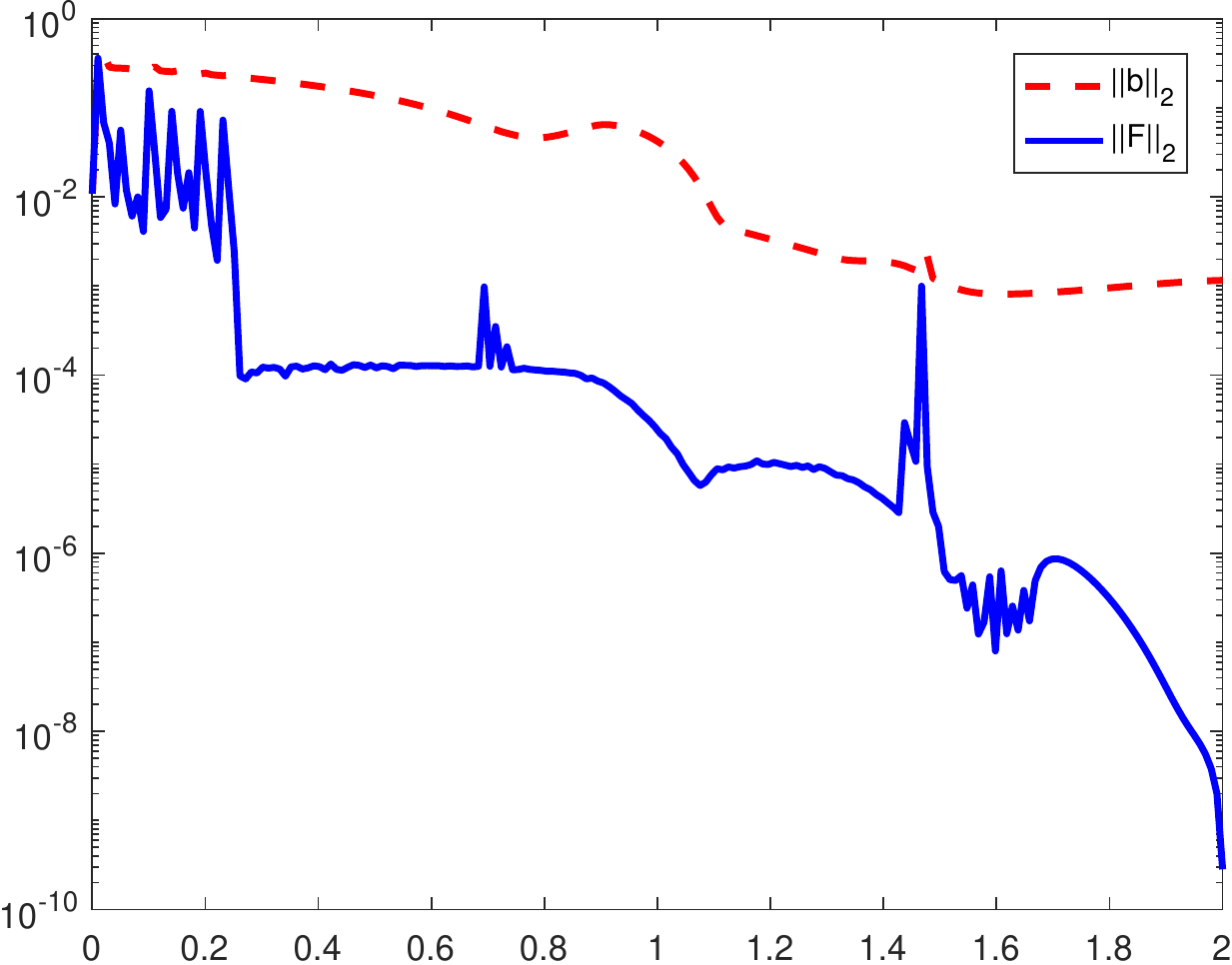}
\caption{Model 2: 2-norm of $F(U_j,x_j,t_j)$ before (dashed line) and
after (solid line) Newton's refinement.}
\label{fig7}
\end{figure}

\begin{figure}
\center
\includegraphics[width=20em]{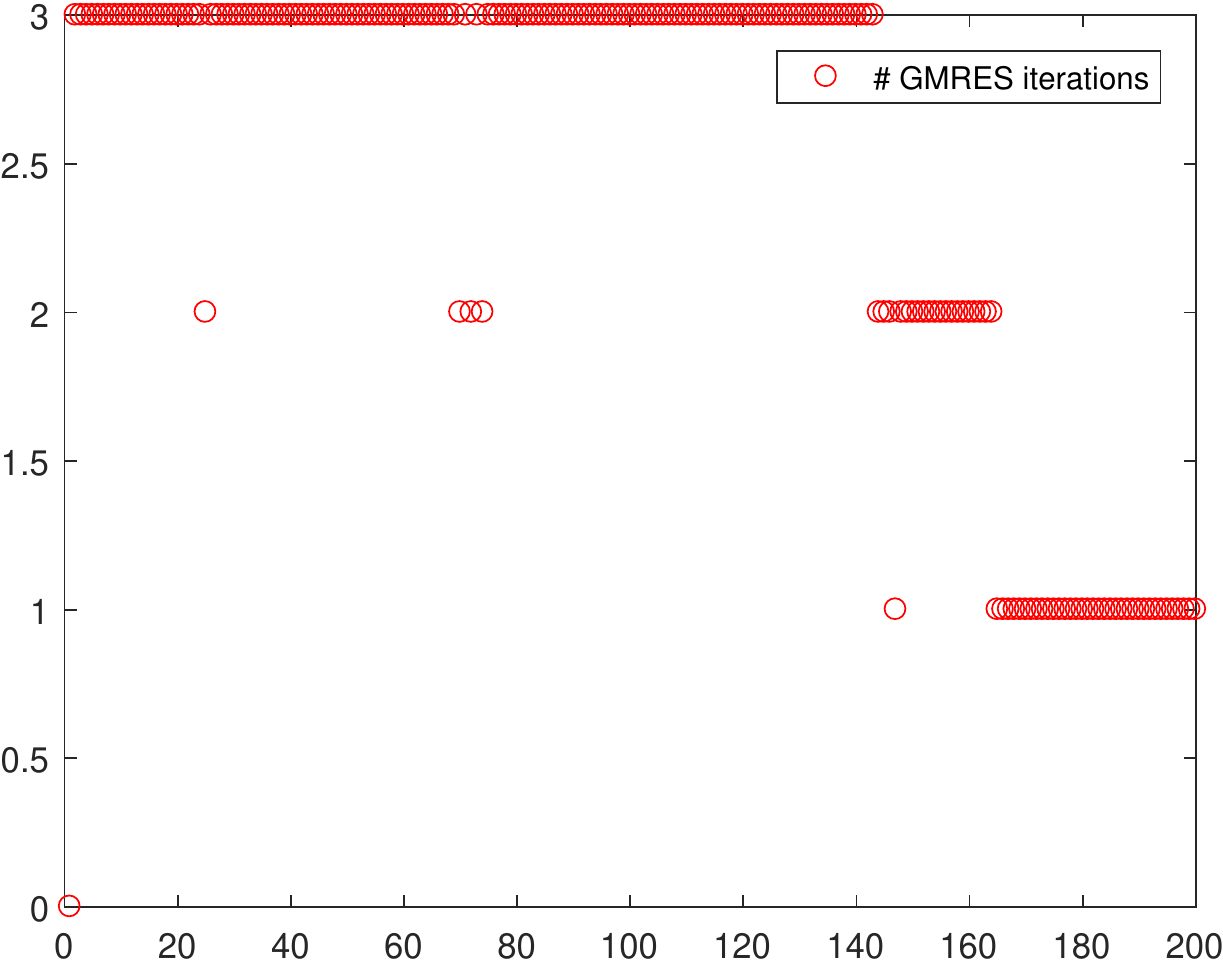}
\caption{Model 2: number of GMRES iterations with preconditioning
(2.5 iterations per step).}
\label{fig8}
\end{figure}

\begin{figure}
\center
\includegraphics[width=20em]{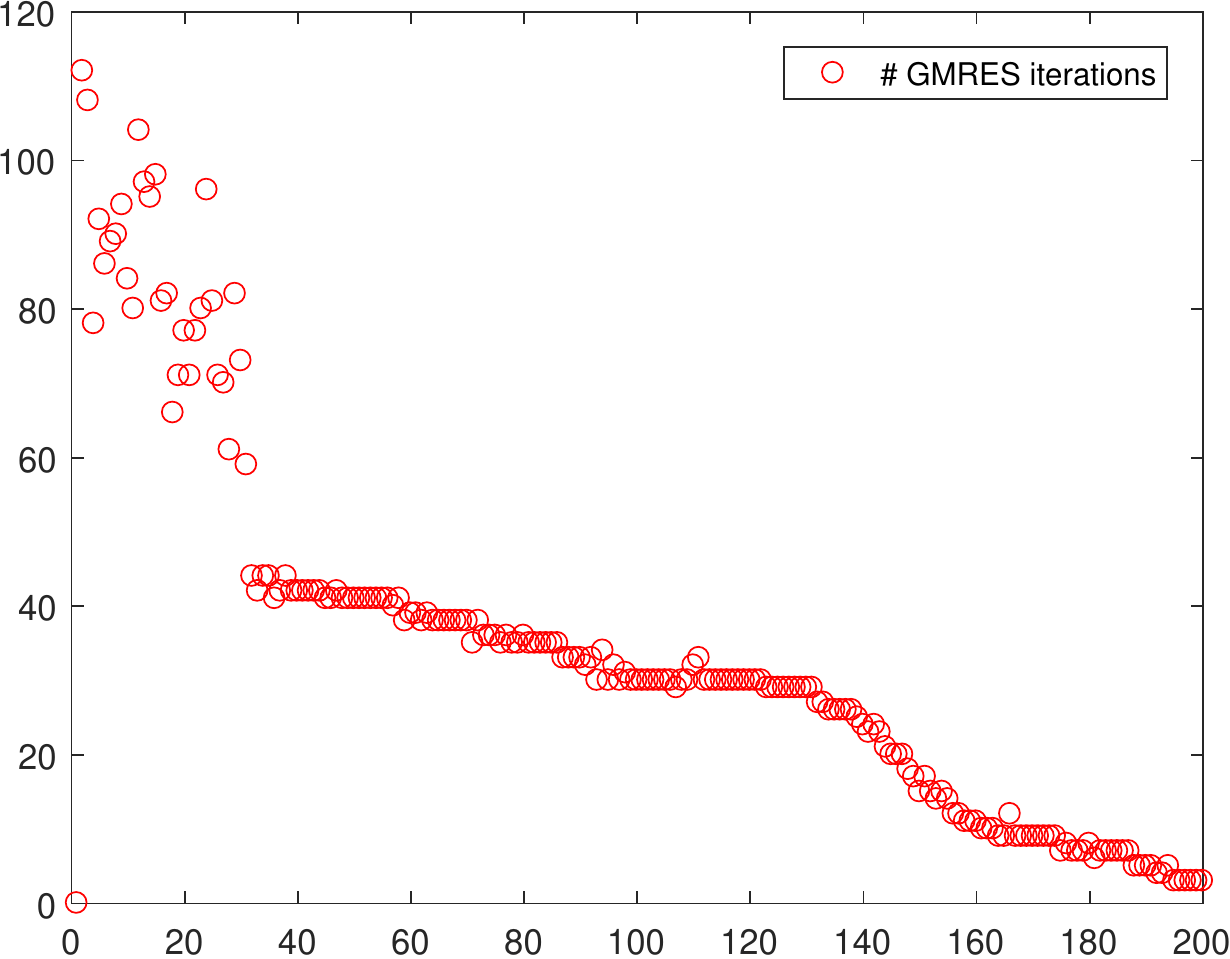}
\caption{Model 2: number of GMRES iterations without preconditioning
(35 iterations per step).}
\label{fig9}
\end{figure}

Figures \ref{fig5}--\ref{fig9} show the numerical results
for Model 2. Again, the computed time to destination equals $t_f=2$, but the system time between
$t_0=0$ and $t_f$ is partitioned into 200 subintervals. The predictive horizon has the length
$N=70$. Our MATLAB program succeeds in computing the control input and trajectory with only
1 Newton's refinement per step by using the shifting along the predictive horizon.
The program fails without the shifting interpolation. The number of GMRES iterations without preconditioning is shown in Fig.~\ref{fig9} and its average value per step is about 35.
The number of GMRES iterations with preconditioning shown in Fig.~\ref{fig8} is about
2.5 per step. Hence the speedup is about 14x.

\begin{figure}
\center
\includegraphics[width=20em]{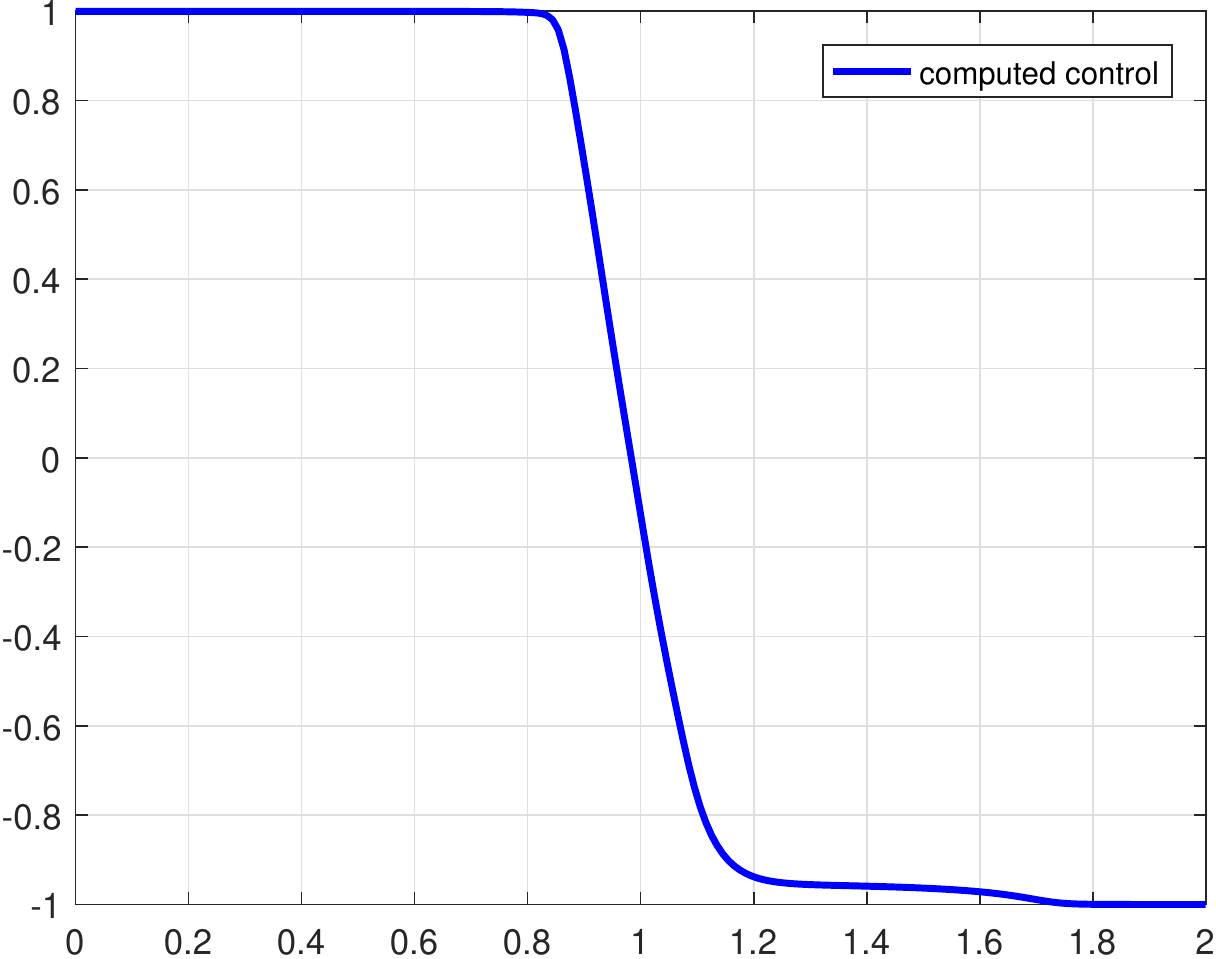}
\caption{Model 2: computed control, $N=20$.}
\label{fig10}
\end{figure}

We have run our program for Model 2, when $N=20$, the time interval $[0,t_f]$ is partitioned
into 200 subintervals, and only 1 Newton's refinement is used. The program fails to compute
without the shifting interpolation and works with the interpolation. Figure~\ref{fig10}
displays the computed control input.

Models 1 and 2 demonstrate different behavior in our numerical experiments.
Model 2 produces the control input with a smoother approximation of the discontinuity than Model 1 in spite of
that its predictive horizon contains 70 grid points against 20 grid points in Model 1.
Model 1 requires application of several Newton's iterations per time step and the shifting
interpolation is not enough for the work without failure. Model 2 works perfectly with
the shifting interpolation and fails with it. The multiple Newton refinements help for Model 2
but the number of required refinements is not small.
Thus, Model 1 with the chosen set of parameters requires multiple Newton iterations, and
Model 2 with its own set of parameters requires the shifting interpolation along the horizon.

\section{Conclusion}

Our numerical experiments demonstrate that the Newton-Krylov method successfully works
even in the cases, when the control input is discontinuous. Both the multiple Newton
refinement and the interpolation by shifting along the predictive horizon can
essentially improve accuracy of the computed solution and avoid failure during
computations. The sparse preconditioner gives very good acceleration, e.g.,
a speedup about 15-30 times. A future research may include more advanced numerical examples
and deeper study of the interpolation by shifting along the predictive horizon.

\paragraph{Acknowledgements.} We are very grateful to Ilya Lashuk for his guidance and help
in using Model 2.

\begin{thebibliography}{99}

\bibitem{AtFa:06}
M. Athans and P. Falb, {\em Optimal control: an introduction to the theory and its
applications}, Dover Publ., 2006.

\bibitem{BrHo:75}
A. E. Bryson, Jr. and Yu-Chi Ho, {\em Applied optimal control: optimization, estimation, and
control}, Taylor \& Francis, 1975.

\bibitem{ChBaLe:12}
D. Chen, L. Bako and S. Lecoeuche, {\em The minimum-time problem for discrete-time linear systems:
a non-smooth optimization approach}, IEEE Int. Conf. Control Appl. (CCA), Dubrovnik, Croatia,
pp.~196--201, 2012.

\bibitem{DiFeHa:09}
M. Diehl, H. J. Ferreau, and N. Haverbeke,
{\em Efficient numerical methods for nonlinear MPC and moving horizon estimation},
in L. Magni et al. (eds.) Nonlinear Model Predictive Control,
LNCIS 384, pp.~391--417, Springer, Heidelberg, 2009.

\bibitem{Fi:88}
A.~F. Filippov, {\em Differential equations with discontinuous
righthand sides}, Springer Netherlands, 1988.

\bibitem{Ga04}
Z. Gao, {\em On discrete time optimal control: a closed-form solution},
in Proc. American Control Conf., pp.~52--58, Boston, MA, 2004.

\bibitem{Gi:13}
P.~Giselsson, {\em Optimal preconditioning and iteration complexity bounds for
gradient-based optimization in model predictive control},
in Proc. American Control Conf., pp.~358--364, Washington D.C., USA, 2013.

\bibitem{Ke:95}
C.~T. Kelly, {\em Iterative methods for linear and nonlinear equations},
SIAM, Philadelphia, PA, 1995.

\bibitem{KnFuMa:15}
A.~Knyazev, Y.~Fujii, and A.~Malyshev, {\em Preconditioned continuation model predictive control},
in Proc. SIAM Conf. Control Appl., pp.~101--108, Paris, France, 2015.

\bibitem{KnMa:15a}
A.~Knyazev and A.~Malyshev,
{\em Preconditioning for continuation model predictive control},
IFAC-PapersOnLine, 48(23) (2015), pp.~191--196.

\bibitem{KnMa:15b}
A.~Knyazev and A.~Malyshev,
{\em Continuation model predictive control on smooth manifolds},
IFAC-PapersOnLine, 48(25) (2015), pp.~126--131.

\bibitem{KnMa:15c}
A.~Knyazev and A.~Malyshev, {\em Efficient particle continuation model predictive control},
IFAC-PapersOnLine, 48(25) (2015), pp.~287--291.

\bibitem{KnMa:16}
A.~Knyazev and A.~Malyshev,
{\em Sparse preconditioning for model predictive control},
in Proc. American Control Conf., pp.~4494--4499, Boston, USA, 2016.

\bibitem{KoFePeDi:15}
D.~Kouzoupis, H.~J. Ferreau, H.~Peyrl, and M.Diehl,
{\em First-order methods in embedded nonlinear model predictive control},
in Proc. European Control Conf., pp.~2622--2627, Linz, Austria, 2015.

\bibitem{Lo:17}
A. Locatelli, {\em Optimal control of a double integrators:
A primer on maximum principle}, Springer Int. Publ., Switzerland, 2017.

\bibitem{Ma:02}
J.~M. Maciejowski, {\em Predictive Control with Constraints},
Prentice-Hall, Englewood Cliffs, NJ, 2002.

\bibitem{Oht:04}
T.~Ohtsuka, {\em A continuation/gmres method for fast computation of nonlinear
  receding horizon control}, Automatica, 40 (2004), pp.~563--574.

\bibitem{Pa:14}
H. Park,
{\em Real-time predictive control of constrained
nonlinear systems using the IPA-SQP approach},
PhD thesis, Univ. Michigan, 2014.

\bibitem{Pon:62}
L.~S. Pontryagin, {\em Ordinary differential equations},
Addison-Wesley, Reading, MA, 1962.

\bibitem{RaWR:98}
C.~V. Rao, S.~J. Wright, and J.~B. Rawlings,
{\em Application of interior-point methods to model predictive control},
J. Optimiz. Theory Appl., 99 (1998), pp.~723--757.

\bibitem{RaMa:13}
J.~B. Rawlings and D.~Q. Mayne,
{\em Model Predictive Control. Theory and Design}, Nob Hill Publishing, Madison, WI, 2013.

\bibitem{Sa:03}
Y. Saad, {\em Iterative Methods for Sparse Linear Systems, 2nd ed.},
SIAM, Philadelphia, 2003.

\bibitem{ShKeCo:10a}
A.~Shahzad, E.~C. Kerrigan, and G.~A. Constantinides,
{\em Preconditioners for inexact interior point methods for predictive
  control}, in Proc. American Control Conf., pp.~5714--5719, Baltimore,
  MD, USA, 2010.

\bibitem{ShKeCo:10b}
A.~Shahzad, E.~C. Kerrigan, and G.~A. Constantinides,
{\em A fast well-conditionend interior point method for predictive
  control}, in Proc. IEEE Conf. Decision Control, pp.~508--513, Atlanta,
  GA, USA, 2010.

\bibitem{ShKeCo:12}
A. Shahzad, E.~C. Kerrigan, and G.~A. Constantinides,
{\em A stable and efficient method for solving a convex quadratic program
  with application to optimal control}, SIAM J. Optim., 22(4) (2012), pp.~1369--1393.

\bibitem{ShOhDi:09}
Y. Shimizu, T. Ohtsuka, and M. Diehl,
{\em A real-time algorithm for nonlinear receding horizon control using
  multiple shooting and continuation/krylov method},
Int. J. Robust Nonlinear Control, 19 (2009), pp.~919--936.

\bibitem{TaOh:04}
T. Tanida and T. Ohtsuka,
{\em Preconditioned c/gmres algorithm for nonlinear receding horizon
  control of hovercrafts connected by a string},
In Proc. IEEE International Conf. Control Appl., pp.~1609--1614,
  Taipei, Taiwan, 2004.

\bibitem{BrDiSW:11}
L. Van den Broeck, M. Diehl, and J. Swevers, {\em A model predictive control
approach for time optimal point-to-point motion control}, Mechatronics,
21(7) (2011), pp.~1203--1212.

\bibitem{WaBo:10}
Y. Wang and S. Boyd,
{\em Fast model predictive control using online optimization},
IEEE Trans. Control Syst. Technology, 18 (2010), pp.~267--278.

\bibitem{ZaBi:09}
V.~M. Zavala and L.~T. Biegler,
{\em Nonlinear programming strategies for state estimation and model
  predictive control}, in L.~Magni et~al., editors, Nonlinear Model Predictive Control,
  LNCIS 384, pp.~419--432, Springer-Verlag, Berlin Heidelberg, 2009.

\end{thebibliography}
\end{document}